\documentclass[11pt]{article}
\usepackage{a4}
\usepackage{amssymb,amsmath,amsthm}
\usepackage{amscd}

\DeclareMathOperator{\st}{st}
\DeclareMathOperator{\lk}{lk}
\DeclareMathOperator{\g}{g}
\DeclareMathOperator{\h}{h}
\DeclareMathOperator{\f}{f}

\DeclareMathOperator{\Ker}{Ker}

\DeclareMathOperator{\chr}{char}
\DeclareMathOperator{\ddim}{dim}
\DeclareMathOperator{\sgn}{sign}

\DeclareMathOperator{\kspan}{span_k}
\DeclareMathOperator{\sspan}{span}
 \DeclareMathOperator{\GIN}{GIN}
\DeclareMathOperator{\ddeg}{deg} \DeclareMathOperator{\mmin}{min}

\DeclareMathOperator{\dom}{dom} \DeclareMathOperator{\res}{res}
\DeclareMathOperator{\ext}{ext}
\DeclareMathOperator{\Stellar}{Stellar}

\theoremstyle{plain}

\newtheorem{theorem}{Theorem}[section]
\newtheorem{proposition}[theorem]{Proposition}
\newtheorem{lemma}[theorem]{Lemma}
\newtheorem*{lemma*}{Auxiliary Lemma}

\newtheorem{corollary}[theorem]{Corollary}

\newtheorem{conjecture}[theorem]{Conjecture}
\newtheorem{problem}[theorem]{Problem}

\theoremstyle{definition}

\theoremstyle{remark}

\newtheorem{remark}[theorem]{Remark}
\newtheorem{remarks}[theorem]{Remarks}

\newtheorem*{acknowledgement}{Acknowledgements}

\title{Lefschetz Properties and Basic Constructions on Simplicial Spheres}
\author{Eric Babson\footnote{Department of Mathematics, UC Davis, Davis USA, E-mail address: babson@math.ucdavis.edu} \ and Eran Nevo\footnote{Department of Mathematics, Cornell
University, Ithaca USA, E-mail address:
eranevo@math.cornell.edu}}

\begin{document}
\maketitle
\begin{abstract}
The well known $g$-conjecture for homology spheres follows from the stronger conjecture that the face ring over the reals of a homology sphere, modulo a linear system of parameters, admits the strong-Lefschetz property. We prove that the strong-Lefschetz property is preserved under the following constructions on homology spheres: join, connected sum, and stellar subdivisions. The last construction is a step towards proving the $g$-conjecture for piecewise-linear spheres.
\end{abstract}

\section{Introduction}\label{sec:Introduction}
Our motivating problem is the following well known $g$-conjecture for spheres, first raised as a question by McMullen for simplicial spheres \cite{McMullen-g-conj}. By \emph{homology sphere}
 we mean a pure simplicial complex
$L$ such that for every face $F\in L$ (including the empty set),
its link $\rm{lk}(F,L):=\{T\in L: T\cap
F=\emptyset,\ T\cup F\in L\}$  has the same homology (say with integer coefficients)
as of a $\rm{dim}(\rm{lk}(F,L))$-sphere. Any simplicial sphere is a homology sphere.
\begin{conjecture}(McMullen \cite{McMullen-g-conj})\label{conj-g}
The $g$-vector of any homology sphere is an $M$-sequence, i.e. is the $f$-vector of a multicomplex.
\end{conjecture}

An algebraic approach to this problem is to associate with a homology sphere $L$ a
standard ring whose Hilbert function is the $g$-vector of
$L$. This was worked out successfully by Stanley \cite{Stanley:NumberFacesSimplicialPolytope-80} in his
celebrated proof of Conjecture \ref{conj-g} for the case where $L$
is the boundary complex of a simplicial polytope. The strong-Lefschetz
theorem for toric varieties associated with rational polytopes,
translates in this case to the following property of face rings,
called \emph{strong-Lefschetz}.

Let $K$ be a $(d-1)$-dimensional simplicial complex on the vertex set $[n]$.
The $i$-th skeleton of $K$ is $K_{i}=\{S\in K:
|S|=i+1\}=K\cap \binom{[n]}{i+1}$, its $f$-\emph{vector} is
$\f(K)=(\f_{-1},\f_0,...,\f_{d-1})$ where $\f_i=|K_i|$,
its $h$-\emph{vector} is
$\h(K)=(\h_0,\h_1,...,\h_{d})$ where
$\h_k=\sum_{0\leq i\leq
k}(-1)^{k-i}\binom{d-i}{k-i}\f_{i-1}$, and in case the $h$-vector is symmetric,  its $g$-\emph{vector} is
$\g(K)=(\g_0,...,\g_{\lfloor d/2\rfloor})$ where
$\g_0=\h_0=1$ and $\g_i=\h_i-\h_{i-1}$ for $1\leq i\leq
\lfloor d/2\rfloor$.

Let $\mathbb{F}$ be a field,
$A=\mathbb{F}[x_1,..,x_n]$ be the polynomial ring over $\mathbb{F}$, where each variable has degree one, and $A_i$ is the degree $i$ part of $A$. The \emph{face ring} of $K$, called also Stanley-Reisner ring, is
$\mathbb{F}[K]=A/I_{K}$ where $I_K$ is the ideal in $A$ generated by
the monomials whose support is not an element of $K$. Let
$\Theta=(\theta_1,..,\theta_d)$ be a linear system of parameters (l.s.o.p. for short) of $\mathbb{F}[K]$ - if $\mathbb{F}$ is infinite
it exists, e.g. \cite[Lemma 5.2]{StanleyGreenBook}, and generic
degree one elements will do. Denote $H(K)=
H(K,\Theta)=\mathbb{F}[K]/(\Theta)=H(K)_0\oplus H(K)_1\oplus...$
where the grading is induced by the degree grading in $A$, and
$(\Theta)$ is the ideal in $\mathbb{F}[K]$ generated by the images
of the elements of $\Theta$ under the projection $A\rightarrow
\mathbb{F}[K]$. $K$ is called \emph{Cohen-Macaulay} (CM for short)
over $\mathbb{F}$ if for an (equivalently, every) l.s.o.p. $\Theta$,
$\mathbb{F}[K]$ is a free $\mathbb{F}[\Theta]$-module. If $K$ is CM
then $\ddim_{\mathbb{F}}H(K)_i=\h_i(K)$. (The converse is also true:
$\h$ is an $M$-vector iff $\h=\h(K)$ for some CM complex $K$
\cite[Theorem 3.3]{StanleyGreenBook}.) For $K$ a CM simplicial
complex with a symmetric $h$-vector, if there exists an l.s.o.p.
$\Theta$ and an element $\omega\in A_1$ such that the multiplication
maps $\omega^{d-2i}: H(K,\Theta)_i\longrightarrow
H(K,\Theta)_{d-i}$, $m\mapsto \omega^{d-2i}m$, are isomorphisms for
every $0\leq i\leq \lfloor d/2\rfloor$, we say that $K$ has the
\emph{strong-Lefschetz property}, or that $K$ is SL (over $\mathbb{F}$).

As was shown by Stanley \cite{St}, for $K$
 the boundary complex of a simplicial rational $d$-polytope $P$,
 the l.s.o.p $\Theta$ induced by the embedding of its vertices in $\mathbb{R}^d$
 and $\omega=\sum_{1\leq i\leq n}x_i$ demonstrate that $K$ is SL over $\mathbb{R}$;
 hence so do generic $(\Theta,\omega)$.

Our main result is that the following constructions on homology spheres preserve the strong-Lefschetz property.

\begin{theorem}\label{thm:connsum,join,Stellar}
Let $K$ and $L$ be homology spheres over a field $\mathbb{F}$, and let $F$ be a face of $K$. Denote by $*$ the join operator, by $\#$ the connected sum operator, and by $\Stellar(F,K)$ the stellar subdivision of $K$ at $F$. The following holds:

(1) If $K$ and $L$ are SL over $\mathbb{F}$ and $\mathbb{F}$ has characteristic zero then $K*L$ is a SL homology sphere (over $\mathbb{F}$).

(2) If $K$ and $L$ have the same dimension and are SL then $K\#L$ is a SL  homology sphere. (True over any field.)

(3) If $K$ and $\lk(F,K)$ are SL over $\mathbb{R}$ then $\Stellar(F,K)$ is a SL homology sphere (over $\mathbb{R}$).    In particular, if $K$ is SL over $\mathbb{R}$ then so is its barycentric subdivision.
\end{theorem}

\begin{remarks}\label{rem:PL-SL}
(1) Replacing the class of homology spheres by the class of piecewise linear (PL) spheres, Theorem \ref{thm:connsum,join,Stellar} still holds. More generally, if $\mathcal{S}$ is a class of simplicial complexes with the SL property, then any complex in its closure w.r.t. join and connected sum is also SL. If $\mathcal{S}$ is closed under links, then any complex in its closure w.r.t. stellar subdivisions is also SL.
\\
(2) Any PL-sphere can be obtained from the boundary of a simplex by a sequence of stellar subdivisions and their inverses (e.g. the survey \cite{Lickorish:SimplicialMoves-1999}). Thus, to prove the $g$-conjecture for PL-spheres it is left to prove that the SL property is preserved under the inverse of stellar subdivisions, in the case of PL-spheres. For arbitrary complexes, the inverse moves may destroy the SL property, which indicates that this direction is more difficult to prove.
\\
(3) A similar result to Theorem \ref{thm:connsum,join,Stellar}(3) was obtained recently, and independently, by Murai \cite{Murai-EdgeDecomposable}, using different ideas: if one assumes that $\lk(F,K)* \partial(F\setminus \{u\})$ is SL for some $u\in F$ instead of that $\lk(F,K)$ is SL, the conclusion $\Stellar(F,K)$ is SL still holds. His proof works for arbitrary field.
Can his proof be used to prove Theorem \ref{thm:connsum,join,Stellar}(3) for arbitrary field?
\\
(4) We use Theorem \ref{thm:connsum,join,Stellar}(1) to prove Theorem \ref{thm:connsum,join,Stellar}(3). Can Murai's result \cite{Murai-EdgeDecomposable} be used to prove the assertion Theorem \ref{thm:connsum,join,Stellar}(1) for arbitrary field?
\end{remarks}

The CM property and the strong-Lefschetz property have equivalent formulations in terms of the combinatorics of the symmetric algebraic shifting of the original simplicial complex \cite{Kalai:Diameter-91} (definitions and further details appear in Section \ref{sec:AlgShifting}).
We consider this reformulation in the context of exterior algebraic shifting, and extend some of our results to this context as well.

This paper is organized as follows:
in Section \ref{sec:SL&Join} we discuss the effect of join on face rings and prove Theorem \ref{thm:connsum,join,Stellar}(1).
In Section \ref{sec:AlgShifting} we give background on algebraic shifting and the interpretation of various Lefschetz properties in terms of shifting.
In Section \ref{sec:SLversusWL} we compare the strong and weak-Lefschetz properties, to be used later in the proof of Theorem \ref{thm:connsum,join,Stellar}(3).
In Section \ref{sec:SL&Stellar} we relate a certain Lefschetz type property, in terms of algebraic shifting (symmetric and exterior), to certain edge contractions, and use it to conclude Theorem \ref{thm:connsum,join,Stellar}(3).
In Section \ref{sec:WL,SL&Connected Sum} we show that connected sum preserves both the strong and weak-Lefschetz properties, also in the exterior algebra context; in particular we prove Theorem \ref{thm:connsum,join,Stellar}(2).


\section{Strong-Lefschetz and join}\label{sec:SL&Join}
The following auxiliary lemma is used in the proof of Theorem \ref{thm:connsum,join,Stellar}(1).
\begin{lemma}\label{lemma:w-invariant subspaces}
Let $K$ be a $(d-1)$-dimensional homology sphere with an
l.s.o.p. $\Theta$ and an SL element $\omega$ over $\mathbb{F}$. Let
$H=\mathbb{F}[K]/(\Theta)$. Then $H$
decomposes into a direct sum of $\mathbb{F}[\omega]$-invariant
spaces, each is of the form $$V_m=\mathbb{F}m\oplus \mathbb{F}\omega
m\oplus...\oplus \mathbb{F}\omega^{d-2i} m$$ for $m\in
\mathbb{F}[K]/(\Theta)$ of degree $i$ for some $0\leq i \leq d/2$.
\end{lemma}
$Proof$: $V_1$ ($1\in H_0$) is an $\mathbb{F}[\omega]$-invariant
space which contain $H_0$. Assume that for $1\leq i\leq d/2$ we have
already constructed a direct sum of $\mathbb{F}[\omega]$-invariant
spaces, $\tilde{V}_{i-1}$, which contains
$\tilde{H}_{i-1}:=H_0\oplus...\oplus H_{i-1}$, in which each $V_m$
contains some nonzero element of $\tilde{H}_{i-1}$. We now extend
the construction to have these properties w.r.t. $\tilde{H}_i$.

Let $W_i:=\ker(\omega^{d-2i+1}:H_i\rightarrow H_{d-i+1})$, and let $m_1,...,m_t$ form a basis (over $\mathbb{F}$) to $W_i$. By definition of $W_i$ each $V_{m_j}$, $1\leq j\leq t$ is $\mathbb{F}[\omega]$-invariant. As $\omega^{d-2i}:H_i\rightarrow H_{d-i}$ is injective, the sum of the $V_{m_j}$'s is direct, denoted by $V_i=\bigoplus_{1\leq j\leq t}V_{m_j}$. Let us check that $V_i\cap \tilde{V}_{i-1}=0$ by showing that its intersection with each $H_l$ is zero. For $l>d-i$ or $l<i$ this is obvious. Otherwise, an element in $V_i\cap \tilde{V}_{i-1}\cap H_l$ is of the form $\omega^{l-i+1}x=\omega^{l-i}y$ where $x\in H_{i-1}$, $y\in W_i$ and $i\leq l\leq d-i$. As $\omega$ is a SL-element, multiplying by $\omega^{d-i+1-l}$, the LHS is nonzero while by definition of $W_i$ the RHS is zero, a contradiction.
We now show that the direct sum in degree $i$ $(V_i\oplus \tilde{V}_{i-1})_i$ equals $H_i$, by computing dimensions:
$\dim_{\mathbb{F}}(\tilde{V}_{i-1})_i=\dim_{\mathbb{F}}(\omega H_{i-1})_i=\h_{i-1}(K)$, and $\dim_{\mathbb{F}}W_i=\h_{i}(K)-\h_{d-i+1}(K)= \h_{i}(K)-\h_{i-1}(K)$ hence $(V_i\oplus \tilde{V}_{i-1})_i=H_i$ and $\tilde{H}_i$ has the desired properties. As the $h$-vector of $K$ is symmetric,  $H=\tilde{H}_{\lfloor d/2\rfloor}$, which completes the proof. $\square$

Recall that the join of two simplicial complexes with disjoint sets of vertices is $K*L:=\{S\cup T: S\in K, T\in L\}$.

\begin{theorem}\label{thm:*}
Let $K$ and $L$ be homology spheres over a field $\mathbb{F}$ on disjoint sets of
vertices, of dimensions $d_K-1,d_L-1$, with l.s.o.p's
$\Theta_K,\Theta_L$ and SL elements $\omega_K,\omega_L$
respectively; over $\mathbb{F}$. Then:

(0) $K*L$ is a homology sphere of dimension $d_K+d_L-1$.

(1) $\Theta_K\biguplus\Theta_L$ is an l.s.o.p for $K*L$ (over
$\mathbb{F}$).

(2) If $\chr(\mathbb{F})=0$ then $\omega_K+\omega_L$ is an SL element of
$\mathbb{F}[K*L]/(\Theta_K\biguplus\Theta_L)$.
\end{theorem}
$Proof:$
(0) is easy and well known; it implies that $K*L$ is CM with a symmetric $h$-vector. We now exhibit a special l.s.o.p. for $K*L$.

For a set $I$ let $A_I:=\mathbb{F}[x_i: i\in I]$ be a polynomial
ring. The isomorphism $A_{K_0}\bigotimes_{\mathbb{F}}A_{L_0}\cong
A_{K_0\biguplus L_0}$, $a_K\otimes a_L\mapsto a_Ka_L$ induces a
structure of an $A=A_{K_0\biguplus L_0}$ module on
$\mathbb{F}[K]\bigotimes_{\mathbb{F}}\mathbb{F}[L]$, isomorphic to
$\mathbb{F}[K*L]$, by $m_K\otimes m_L\mapsto m_Km_L$ and
$(a_K\otimes a_L)(m_K\otimes m_L)=a_Km_K\otimes a_Lm_L$. (E.g.
$a_K\in A_{K_0}\subseteq A$ acts like $a_K\otimes 1$ on
$\mathbb{F}[K]\bigotimes_{\mathbb{F}}\mathbb{F}[L]$. )

The above isomorphism induces an isomorphism of $A$-modules
\begin{equation}\label{eq:JoinIsom}
\mathbb{F}[K*L]/(\Theta_K\biguplus\Theta_L)\mathbb{F}[K*L] \cong
\mathbb{F}[K]/(\Theta_K)\mathbb{F}[K] \bigotimes_{\mathbb{F}}
\mathbb{F}[L]/(\Theta_L)\mathbb{F}[L],
\end{equation}
proving (1).

By Lemma \ref{lemma:w-invariant subspaces},
$\mathbb{F}[K]/(\Theta_K)$ decomposes into a direct sum of
$\mathbb{F}[\omega_K]$-invariant spaces, each is of the form
$V_m=\mathbb{F}m\bigoplus \mathbb{F}\omega_K m\bigoplus...\bigoplus
\mathbb{F}\omega_K^{d_K-2i} m$ for $m\in \mathbb{F}[K]/(\Theta_K)$
of degree $i$ for some $0\leq i \leq d_K/2$; and similarly for
$\mathbb{F}[L]/(\Theta_L)$.

First let us consider the case $\mathbb{F}=\mathbb{R}$:
the $\mathbb{R}[\omega_K]$-module $V_m$ is isomorphic to the
$\mathbb{R}[\omega]$-module $\mathbb{R}[\partial
\sigma^{d_K-2i}]/(\theta)$ by $\omega_K\mapsto \omega$ and $m\mapsto
1$, where $\sigma^j$ is the $j$-simplex, $\theta$ is an l.s.o.p.
induced by the positions of the vertices in an embedding of
$\sigma^{d_K-2i}$ as a full dimensional geometric simplex in
$\mathbb{R}^{d_K-2i}$ with the origin in its interior, and
$\omega=\sum_{v\in \sigma_0}x_v$ is an SL element for
$\mathbb{R}[\partial \sigma^{d_K-2i}]/(\theta)$. Thus, to prove (2) for $\mathbb{F}=\mathbb{R}$
it is enough to prove it for the join of boundaries of two simplices
with l.s.o.p.'s as above and the SL elements having weight $1$ on
each vertex of the ground set.

Note that the join $\partial \sigma^{k}*\partial \sigma^{l}$ is
combinatorially isomorphic to the boundary of the polytope
$P:=\rm{conv}(\sigma^{k}\cup_{\{0\}} \sigma^{l})$ where $\sigma^{k}$
and $\sigma^{l}$ are embedded in orthogonal spaces and intersect
only in the origin which is in the relative interior of both.
McMullen's proof of the $g$-theorem for simplicial polytopes
\cite{McMullen-g-proof1,McMullen-g-proof2} states that $\sum_{v\in
P_0}x_v = \omega_{\partial \sigma^{k}}+\omega_{\partial \sigma^{l}}$
is indeed an SL element of $\mathbb{R}[\partial \sigma^{k}*\partial
\sigma^{l}]/(\Theta_{\partial P})$ where $\Theta_{\partial P}$ is
the l.s.o.p. induced by the positions of the vertices in the
polytope $P$. By the definition of $P$, $\Theta_{\partial
P}=\Theta_{\partial \sigma^{k}}\uplus \Theta_{\partial \sigma^{l}}$.
Thus (2) is proved for $\mathbb{F}=\mathbb{R}$.

For a general field with $\chr(\mathbb{F})=0$, notice that
$V_{m}$ as above ($m\in \mathbb{F}[K]$ homogenous) is isomorphic as an $\mathbb{F}[\omega_K]$-module to $\mathbb{F}[\omega_K]/\omega_K^{d_K-2i+1}\mathbb{F}[\omega_K]$, hence for $\omega=\omega_K+\omega_L$ we get an isomorphism of $\mathbb{F}[\omega]$ modules
$$\mathbb{F}[V_{m(K)}*V_{m(L)}]\cong \frac{\mathbb{F}[\omega_K]}{(\omega_K^{d_K-2i_K+1})\mathbb{F}[\omega_K]}\bigotimes_{\mathbb{F}} \frac{\mathbb{F}[\omega_L]}{(\omega_L^{d_L-2i_L+1})\mathbb{F}[\omega_L]}.$$
Picking the basis $\{\omega_K^l\otimes\omega_L^j:\ 0\leq l\leq d_K-2i_K,0\leq j\leq d_L-2i_L\}$ for the module on the RHS, we see that the representing matrix of the map $\omega^{d_K+d_L-2i}:(\mathbb{F}[V_{m(K)}*V_{m(L)}])_i\rightarrow (\mathbb{F}[V_{m(K)}*V_{m(L)}])_{d_K+d_L-i}$ consist of integer entries (all entries are binomials). The case $\mathbb{F}=\mathbb{R}$ shows that its determinant is nonzero, hence (2) follows for every field of characteristic zero.
$\square$

In particular, Theorem \ref{thm:*} implies Theorem \ref{thm:connsum,join,Stellar}(1). Similarly, as the join of PL spheres is a PL sphere, Remark \ref{rem:PL-SL}(1) follows in the same manner.
$\square$

\begin{remarks}
(1) As a nonzero multiple of an SL element is again SL,
then in Theorem \ref{thm:*}(2) any element $a\omega_K+b\omega_L$
where $a,b\in \mathbb{F}$, $ab\neq 0$, will do.

(2) A closer look at the integer matrix used in the proof shows that if $\chr(\mathbb{F})\neq 0$ then there exist simplices $\sigma^{d_K},\sigma^{d_L}$ such that for any l.s.o.p's $\Theta_K,\Theta_L$ of the face rings of their boundaries, respectively, there is no SL-element for $\mathbb{F}[\partial \sigma^{d_K}* \partial \sigma^{d_L}]/(\Theta_K\cup \Theta_L)$.  On the other hand, for strongly edge decomposable complexes, introduced in \cite{Nevo-VK}, Murai proved recently, see \cite[Corollary 3.5]{Murai-EdgeDecomposable}, that the SL property holds over any field. The join of boundaries of two simplices is strongly edge decomposable (identify a pair of vertices, one from each simplex, to obtain the boundary of a simplex), hence for some other l.s.o.p $\Theta$, $\mathbb{F}[\partial \sigma^{d_K}* \partial \sigma^{d_L}]/(\Theta)$ has an SL-element. This raises the following question:
\end{remarks}
\begin{problem}
Does Theorem \ref{thm:connsum,join,Stellar}(1) hold for a field of arbitrary characteristic?
Can the results in \cite{Murai-EdgeDecomposable} be used to prove this?
\end{problem}

\section{Algebraic shifting}\label{sec:AlgShifting}
Let $<$ denote the usual order on the natural numbers. A simplicial
complex $K$ with vertices $[n]=\{1,2,...,n\}$ is \emph{shifted} if for every
$i<j$ and $j\in S\in K$, also $(S\setminus\{j\})\cup \{i\}\in K$.

Algebraic shifting is an operator associating with each simplicial complex a shifted simplicial complex.
It has two versions - exterior and symmetric, both introduced by Kalai. Various invariants of the original complex, like its $f$-vector and Betti numbers, can be read off from its shifting. For a survey on algebraic shifting see Kalai \cite{skira}. For completeness we give now the definitions of exterior and symmetric shifting.
\\
\textbf{Exterior shifting.}
Let $\mathbb{F}$ be a field and let $k$ be a field extension of
$\mathbb{F}$ of transcendental degree $\geq n^2$ (e.g.
$\mathbb{F}=\mathbb{Q}$ and $k=\mathbb{R}$, or
$\mathbb{F}=\mathbb{Z}_2$ and $k=\mathbb{Z}_2(x_{ij})_{1\leq i,j\leq
n}$ where $x_{ij}$ are intermediates). Let $V$ be an $n$-dimensional
vector space over $k$ with basis $\{e_{1},\dots,e_{n}\}$. Let
$\bigwedge V$ be the graded exterior algebra over $V$. Denote
$e_{S}=e_{s_{1}}\wedge\dots\wedge e_{s_{j}}$ where $S=
\{s_{1}<\dots<s_{j}\}$. Then $\{e_{S}: S\in (_{\ j}^{[n]})\}$ is a
basis for $\bigwedge^{j} V$. Note that as $K$ is a simplicial
complex, the ideal $(e_{S}:S\notin K)$ of $\bigwedge V$ and the
vector subspace $\rm{span}\{e_{S}:S\notin K\}$ of $\bigwedge V$
consist of the same set of elements in $\bigwedge V$. Define the
exterior algebra of $K$ by
$$\bigwedge K=(\bigwedge V)/(e_{S}:S\notin K).$$ Let
$\{f_{1},\dots,f_{n}\}$ be a basis of $V$, generic over $\mathbb{F}$
with respect to $\{e_{1},\dots,e_{n}\}$, which means that the
entries of the corresponding transition matrix $A$ ($e_{i}A=f_{i}$
for all $i$) are algebraically independent over $\mathbb{F}$. Let
$\tilde{f}_{S}$ be the image of $f_{S}\in \bigwedge V$ in
$\bigwedge K$. Let $<_L$ be the lexicographic order on equal sized
subsets of $\mathbb{N}$, i.e. $S<_LT$ iff $\mmin(S\triangle T)\in
S$. Define
$$\Delta^e(K)=\Delta^e_{A}(K)=\{S: \tilde{f}_{S}\notin \rm{span}\{\tilde{f}_{S'}:S'<_{L}S\}\}$$
to be the exterior shifting of $K$, introduced by Kalai \cite{55}.
 The construction is canonical, i.e. it is independent of the choice of the generic matrix
 $A$, and for a permutation $\pi:[n]\rightarrow [n]$ the induced simplicial complex
  $\pi(K)$ satisfies $\Delta^e(\pi(K))=\Delta^e(K)$.
 It results in a shifted simplicial complex,
 having the same face vector and Betti vector as $K$ \cite{BK}.
\\
\textbf{Symmetric shifting.}
let us look on the face ring
(Stanley-Reisner ring) of $K$ $k[K]=k[x_{1},..,x_{n}]/I_{K}$ where
$I_{K}$ is the homogenous ideal generated by the monomials whose
support is not in $K$, $\{\prod_{i\in S}x_i:\ S\notin K\}$. $k[K]$
is graded by degree. Let $\mathbb{F}\subseteq k$ be fields as before
and let $y_{1},\dots,y_{n}$ be generic linear combinations of
$x_{1},\dots,x_{n}$ w.r.t. $\mathbb{F}$. We choose a basis for each
graded component of $k[K]$, up to degree $\ddim(K)+1$, from the
canonic projection of the monomials in the $y_{i}$'s on $k[K]$, in
the greedy way:
$$\GIN (K)=\{m: \tilde{m}\notin \kspan   \{\tilde{m'}:\ddeg (m')=\ddeg (m), m'<_{L}m\}\}$$
where $\prod y_{i}^{a_{i}}<_{L}\prod y_{i}^{b_{i}}$ iff for $j=\mmin
\{i: a_{i}\neq b_{i}\}$ $a_{j}>b_{j}$. The combinatorial
information in $\GIN (K)$ is redundant: if $m\in \GIN (K)$ is of
degree $i\leq \ddim (K)$ then $y_{1}m,..,y_{i}m$ are also in $\GIN
(K)$. Thus, $\GIN (K)$ can be reconstructed from its monomials of
the form $m=y_{i_{1}}\cdot y_{i_{2}}\cdot..\cdot y_{i_{r}}$ where
$r\leq i_{1}\leq i_{2}\leq..\leq i_{r}$, $r\leq \ddim(K)+1$. Denote
this set by $\rm{gin} (K)$, and define
$S(m)=\{i_{1}-r+1,i_{2}-r+2,..,i_{r}\}$ for such $m$. The collection
of sets
$$\Delta^{s}(K)=\cup \{S(m): m\in \rm{gin}(K)\}$$
carries the same combinatorial information as $\GIN (K)$. $\Delta^{s}(K)$ is a
simplicial complex. Again, the construction is canonic, in
the same sense as for exterior shifting.
If $k$ has characteristic zero then $\Delta^{s}(K)$ is shifted \cite{Kalai-SymmMatroids}.
\\
\textbf{Lefschetz properties via shifting.}
$K$ is CM (over $\mathbb{F}$) iff $\Delta^s(K)$ is pure (i.e. all its maximal faces have the same size) and the following condition holds
\begin{equation}\label{eq:CM}
S\in \Delta^s(K), |S|=k \ \Rightarrow \ [d-k]\cup S\in \Delta^s(K).
\end{equation}
To see this take the first $d$ elements in a generic basis $\{y_1,...,y_d\}$ to be an l.s.o.p. for $K$.

Further, let
$\Delta(d,n)$ be the pure $(d-1)$-dimensional simplicial complex
with set of vertices $[n]$ and facets $\{S: S\subseteq [n], |S|=d,\
k\notin S\Rightarrow [k+1,d-k+2]\subseteq S\}$. Equivalently,
$\Delta(d,n)$ is the maximal pure $(d-1)$-dimensional simplicial
complex with vertex set $[n]$ which does not contain any of the sets
$T_d,...,T_{\lceil d/2\rceil}$, where
\begin{equation}\label{eq:T}
T_{d-k}=\{k+2,k+3,...,d-k,d-k+2,d-k+3,...,d+2\},\  0\leq k\leq \lfloor d/2\rfloor.
\end{equation}
Note that $\Delta(d,n)\subseteq\Delta(d,n+1)$, and define $\Delta(d)=\cup_n\Delta(d,n)$.
For $K$ a CM
$(d-1)$-dimensional complex with symmetric $h$-vector,
$\Delta^s(K)\subseteq \Delta(d)$ is equivalent to $K$ being SL.
To see this, take the $(d+1)$'th element in a generic basis, $y_{d+1}$, to be the strong-Lefschetz element: indeed, $\Delta^s(K)\subseteq \Delta(d)$ iff non of the monomials
$y_{d+1}^{d-2k-1}y_{d+2}^{k+1}$ are in $GIN(K)$ (where $k=0,1,...$), iff
the maps $y_{d+1}^{d-2k}: H(K)_k\longrightarrow H(K)_{d-k}$ are onto for $0\leq k\leq \lfloor d/2\rfloor$, and
when $h(K)$ is symmetric this happens iff these maps are
isomorphisms.

Let $\Delta(K)$ refer to both symmetric and exterior shifting.
Kalai refers to the relation
\begin{equation}\label{eq:shUBT}
\Delta(K)\subseteq \Delta(d)
\end{equation}
as the \emph{shifting theoretic upper bound theorem}. To justify the name, note that the
boundary complex of the cyclic $d$-polytope on $n$ vertices, denoted
by $C(d,n)$, satisfies $\Delta^s(C(d,n))=\Delta(d,n)$. This follows
from the fact that $C(d,n)$ is SL. Recently Murai \cite{MuraiCyclic}
proved that also $\Delta^e(C(d,n))=\Delta(d,n)$, as was conjectured
by Kalai \cite{skira}. It follows that if $K$ has $n$ vertices and
(\ref{eq:shUBT}) holds, then the $f$-vectors satisfy $f(K)\leq
f(C(d,n))$ componentwise.

For $K$ as above (CM with symmetric $h$-vector), weaker than the strong-Lefschetz property is to
require only that multiplications $y_{d+1}:
H(K)_{i-1}\longrightarrow H(K)_{i}$ are injective for $1\leq i\leq
\lceil d/2\rceil$ and surjective for $\lceil d/2\rceil<i\leq d$,
usually called in the literature the \emph{weak-Lefschetz} property (WL for short). Even weaker is just to
require that multiplications $y_{d+1}: H(K)_{i-1}\longrightarrow
H(K)_{i}$ are injective for $1\leq i\leq \lfloor d/2\rfloor$, called here WWL property. (Injectivity for $i\leq \lceil d/2\rceil$ in the case of
homology spheres implies also surjective maps for $\lceil
d/2\rceil<i\leq d$ as was noticed by Swartz; see the proof of Theorem \ref{thmSwartz} below.)
The WWL property is equivalent to the following, in the case of symmetric
shifting \cite{BK-homology}:
\begin{eqnarray}\label{eq:shWL}
\ S\in \Delta(K), |S|=k \ \Rightarrow \ [d-k]\cup S\in \Delta(K), \nonumber
\\
\ S\in \Delta(K), |S|=k<\lfloor d/2\rfloor \ \Rightarrow \ \{d-k+1\}\cup S\in \Delta(K).
\end{eqnarray}
The first condition holds when $K$ is CM, and the second condition holds iff $K$
is WWL. As was noticed in \cite{BK-homology}, (\ref{eq:shWL}) is
implied by requiring that $\Delta(K)$ is pure and every
$S\in\Delta(K)$ of size less than $\lfloor d/2\rfloor$ is contained
in at least $2$ facets of $\Delta(K)$.

Note that if $L$ is a homology sphere, it is in particular CM with a
symmetric $h$-vector. If in addition it is WWL, then in the standard ring
$S(L)=\mathbb{F}[L]/(y_1,...,y_{d+1})=H(L,\{y_1,...,y_d\})/(y_{d+1})=S_0\oplus S_1\oplus...$ the following holds:
$g_i(L)=\ddim_\mathbb{F}S_i$ for all $0\leq i\leq \lfloor
d/2\rfloor$, and Conjecture \ref{conj-g} holds for $L$.

We summarize the discussion above in the following hierarchy of conjectures, where assertion $(i)$ implies assertion $(i+1)$:
\begin{conjecture}\label{conj-g-hierarchy}
Let $L$ be a homology $(d-1)$-sphere. Then:

(1) If $S\in\Delta(L)$, $|S|=k\leq \lfloor d/2\rfloor$ and $S\cap[d-k+1]=\emptyset$ then $S\cup[k+2,d-k+1]\in\Delta(L)$.

This is equivalent to $\Delta(K)\subseteq \Delta(d)$, and in the
symmetric case this is equivalent to $L$ being SL.

(2) If $S\in\Delta(L)$, $|S|=k< \lfloor d/2\rfloor$ and
$S\cap[d-k+1]=\emptyset$ then $S\cup[\lceil
d/2\rceil+2,d-k+1]\in\Delta(L)$. In the symmetric case this is
equivalent to $L$ being WWL.

(3) $g(L)$ is an $M$-vector.
\end{conjecture}

\section{Strong Lefschetz versus weak-Lefschetz}\label{sec:SLversusWL}
Examples of Gorenstein algebras admitting the weak-Lefschetz property but not the strong-Lefschetz property were found in \cite[Example 4.3]{HarimaMiglioreNagelWatanabe}. For Gorenstein algebras arising as face rings of homology spheres the SL property is conjectured to hold. Does it follow from the (conjectured) WL property for homology spheres? We end this section with a result in this direction, to be used later in the proof of Theorem \ref{thm:connsum,join,Stellar}(3).

Consider the multiplication maps $\omega_i:
H(K,\Theta)_{i}\longrightarrow H(K,\Theta)_{i+1}$, $m\mapsto
\omega_i m$ where $\omega_i \in A_1$. Let $\ddim(K)=d-1$. Denote by
$\Omega_{WL}(K,i)$ the set of all $(\Theta,\omega_i)\in
A_1^{\ddim(K)+2}$ such that $\Theta$ is an l.s.o.p. of
$\mathbb{F}[K]$, $\mathbb{F}[K]$ is a free
$\mathbb{F}[\Theta]$-module, and $\omega_i: H(K)_{i}\longrightarrow
H(K)_{i+1}$ is injective for $i< d/2$ and surjective for $i\geq
d/2$. Denote by $\Omega_{SL}(K,i)$ the set of all
$(\Theta,\omega)\in A_1^{d+1}$ such that $\Theta$ is an
l.s.o.p. of $\mathbb{F}[K]$, $\mathbb{F}[K]$ is a free
$\mathbb{F}[\Theta]$-module, and $\omega^{d-2i}:
H(K)_{i}\longrightarrow H(K)_{d-i}$ is injective ($0\leq i\leq
\lfloor d/2\rfloor$).
If $\Omega_{SL}(K,i)\neq \emptyset$ we say that $K$ is $i$-Lefschetz and for $(\Theta,\omega)\in \Omega_{SL}(K,i)$ that $H(K,\Theta)$ is $i$-Lefschetz with an $i$-Lefschetz element $\omega$.
For $d$ odd $\Omega_{WL}(K,\lfloor
d/2\rfloor)=\Omega_{SL}(K,\lfloor d/2\rfloor)$, which we simply
denote by $\Omega(K,\lfloor d/2\rfloor)$.

The following is well known, see e.g. \cite[Proposition 3.6]{Swartz}
for the case $\Omega_{SL}(K,i)$; similar arguments can be used to prove the
same conclusion for $\Omega_{WL}(K,i)$.
\begin{lemma}\label{lemOmegaZariski}
For every simplicial complex $K$ and for every $i$,
$\Omega_{WL}(K,i)$ is a Zariski open set. For $0\leq i\leq \lfloor
\frac{\ddim(K)+1}{2}\rfloor$, $\Omega_{SL}(K,i)$ is a Zariski open
set. (They may be empty, e.g. if $K$ is not pure.)
\end{lemma}
\begin{theorem}(Swartz)\label{thmSwartz}
Let $d\geq 1$. If for every homology $2d$-sphere $L$, $\Omega(L,d)$
is nonempty, then for every $t>2d$ and for every homology $t$-sphere
$K$, $\Omega_{WL}(K,m)$ is nonempty for every $m\leq d$.
In particular, the condition implies the WL property for homology spheres, hence Conjecture \ref{conj-g} would follow.
\end{theorem}
$Proof$: By \cite[Theorem 4.26]{Swartz-SpheresToManifolds} and
induction on $t$, $\Omega_{WL}(K,(t+1)-(d+1))$ is nonempty, i.e.
multiplication $\omega :H(K)_{t-d}\rightarrow H(K)_{t-d+1}$ is
surjective for a generic l.s.o.p. and $\omega\in A_1$. As the ring $H(K)$ is
standard, $\Omega_{WL}(K,(t+1)-(m+1))$ is nonempty for every $m\leq
d$. Hence, for the canonical module $\Omega(K)$, multiplication by a
generic degree $1$ element
$\omega:(\Omega(K)/\Theta\Omega(K))_{m}\rightarrow
(\Omega(K)/\Theta\Omega(K))_{m+1}$ is injective in the first $d$
degrees. As $K$ is a homology sphere, $\Omega(K)\cong
\mathbb{R}[K]$ as graded $A$-modules up to a shift in grading (e.g. \cite{StanleyGreenBook}), hence $\Omega_{WL}(K,m)$ is nonempty for every
$m\leq d$.
Combined with Lemma \ref{lemOmegaZariski}, and the fact that a
finite intersection of Zariski nonempty open sets is nonempty, if
the conditions of Theorem \ref{thmSwartz} are met for every $d\geq
1$ then every homology sphere is WL, and hence Conjecture
\ref{conj-g} follows.
$\square$

We wish to show further, that if all even dimensional homology spheres
satisfy the condition in Theorem \ref{thmSwartz} then all homology spheres
are SL. The following result aims at this direction. If one extends its conclusion for \emph{every} l.s.o.p. of $S*\partial
\sigma$, then indeed WL would imply SL for homology spheres.

\begin{lemma}\label{lem:WL->SL}
Let $S$ be a homology sphere with an l.s.o.p.
$\Theta_S$ over a field $\mathbb{F}$ of characteristic zero.
If $H(S,\Theta_S)$ is $(\lfloor
\frac{\ddim S +1}{2}\rfloor)$-Lefschetz but not SL then there exists
a simplex $\sigma$ such that the homology sphere $S*\partial\sigma$ is of even dimension
$2j$, and for every l.s.o.p. $\Theta_{\partial \sigma}$ of $\partial
\sigma$, $\mathbb{F}[S*\partial \sigma]/(\Theta_S \cup
\Theta_{\partial \sigma})$ has no $j$-Lefschetz element; in
particular $S*\partial\sigma$ is not WL.
\end{lemma}

$Proof:$ Denote the dimension of $S$ by $d-1$ and recall that
$A_{S_0}=\mathbb{F}[x_v:v\in S_0]$. By Lemma \ref{lemOmegaZariski}
$\Omega_{SL}(S,i)$ is a Zariski open set for every $0\leq i\leq
\lfloor d/2\rfloor$. The assumption that $S$ is not SL (but is
$(\lfloor \frac{d}{2}\rfloor)$-Lefschetz) implies that there exists
$0\leq i_0\leq \lfloor d/2\rfloor -1$ such that
$\Omega_{SL}(S,i_0)=\emptyset$ (as a finite intersection of Zariski
nonempty open sets is nonempty). Hence, for the fixed l.s.o.p.
$\Theta_S$ and every $\omega_S \in (A_{S_0})_1$, there exists $0\neq
m=m(\omega_S)\in H_{i_0}(S)$ such that $\omega_S^{d-2i_0}m=0$.

Let $T=S*\partial \sigma$ where $\sigma$ is the $(d-2i_0
-1)$-simplex. Note that $\ddim(\sigma)\geq 1$, hence $\partial \sigma\neq
\emptyset$. Then $T$ is a homology sphere of even dimension $2d-2i_0-2$. We have seen
(Theorem \ref{thm:*}) that for any l.s.o.p. $\Theta_{\partial
\sigma}$ of $\partial \sigma$, $\Theta_T:=\Theta_S \cup
\Theta_{\partial \sigma}$ is an l.s.o.p. of $T$. Every $\omega_T\in
(A_{T_0})_1$ has a unique expansion
$\omega_T=\omega_S+\omega_{\partial \sigma}$ where $\omega_S\in
(A_{S_0})_1$ and $\omega_{\partial \sigma} \in (A_{\partial
\sigma_0})_1$. Recall the isomorphism (\ref{eq:JoinIsom}) of
$A_{T_0}$-modules $\mathbb{F}[T]/(\Theta_T)\cong
\mathbb{F}[S]/(\Theta_S)\otimes_{\mathbb{F}} \mathbb{F}[\partial
\sigma]/(\Theta_{\partial \sigma})$. Let $m(\omega_T)\in
(\frac{\mathbb{F}[T]}{(\Theta_T)})_{d-i_0-1}$ be
$$m(\omega_T):=\sum_{0\leq j\leq d-2i_0-1}(-1)^j \omega_S^{d-2i_0-1-j}m\otimes \omega_{\partial \sigma}^j 1.$$
Note that the sum $\omega_T m(\omega_T)$ is telescopic, thus
$\omega_T m(\omega_T)=\omega_S^{d-2i_0}m\otimes 1+ (-1)^{d-2i_0-1}
m\otimes \omega_{\partial \sigma}^{d-2i_0}1= 0+0=0$. For a generic
$\omega_T$, the projection of $\omega_{\partial \sigma}$ on
$\mathbb{F}[\partial \sigma]/(\Theta_{\partial \sigma})$ is nonzero,
hence so is the projection of $\omega_{\partial \sigma}^{d-2i_0-1}$,
and we get that $m(\omega_T)\neq 0$. Thus, Zariski topology tells us
that for \emph{every} $\omega_T\in (A_{T_0})_1$, there exists $0\neq
m(\omega_T)\in (\frac{\mathbb{F}[T]}{(\Theta_T)})_{d-i_0-1}$ such
that $\omega_T m(\omega_T)=0$. $\square$


\section{Lefschetz properties and Stellar subdivisions}\label{sec:SL&Stellar}

Roughly speaking, we will show that Stellar subdivisions preserve
the SL property.
\begin{proposition}\label{propAlgContruction}
Let $K$ be a simplicial complex. Let $K'$ be obtained from $K$ by
identifying two distinct vertices $u$ and $v$ in $K$, i.e. $K'=\{T:
u\notin T \in K\}\cup\{(T\setminus \{u\})\cup \{v\}: u\in T \in
K\}$. Let $d\geq 2$. Assume that $\{d+2,d+3,...,2d+1\}\notin
\Delta(K')$ and that $\{d+1,d+2,...,2d-1\}\notin \Delta(\lk(u,K)\cap
\lk(v,K))$. Then $\{d+2,d+3,...,2d+1\}\notin \Delta(K)$. (Shifting
is over $\mathbb{R}$.)
\end{proposition}
The case $d=2$ and $\ddim(K)=1$ of this proposition was proved by Whiteley \cite{Wh} in the symmetric case. The relation between symmetric shifting and rigidity of graphs, discussed in Lee~\cite{Lee}, is used to translate his result to algebraic shifting terms.
\\
\\
\emph{Proof for symmetric shifting}:
Let $\psi: K_0\longrightarrow \mathbb{R}^{2d}$ be a generic
map, i.e. all minors of the representing matrix w.r.t. a fixed
basis are nonzero. It induces the following map:
\begin{eqnarray}\label{eq:HighRigidity}
\psi_K^{2d}: \oplus_{T\in K_{d-1}}\mathbb{R}T \longrightarrow
\oplus_{F\in \binom{K_0}{d-1}}\mathbb{R}^{2d}/\sspan(\psi(F)),\nonumber \\
1T\mapsto \sum_{F\in \binom{K_0}{d-1}}\delta_{F\subseteq T} \overline{\psi(T\setminus F)} F
\end{eqnarray}
where $\delta_{F\subseteq T}$ equals $1$ if $F\subseteq T$ and $0$
otherwise.

Recall that $\{d+2,d+3,...,2d+1\}\notin \Delta^s(K)$ iff $y_{2d+1}^d
\notin GIN(K)$, where $Y=\{y_i\}_i$ is a generic basis for $A_1$,
$A=\mathbb{R}[x_v: v\in K_0]$.
 By Lee \cite[Theorems 10,12,15]{Lee} and Tay, White and
Whiteley \cite[Proposition 5.2]{TayWhiteWhiteley-skel1}, $y_{2d+1}^d
\notin GIN(K)$ iff $\Ker \phi_K^{2d}=0$ for some $\phi:
K_0\longrightarrow \mathbb{R}^{2d}$ (equivalently, every $\phi$ in
some Zariski non-empty open set of maps).

Consider the following degenerating map: for $0<t\leq 1$ let
$\psi_t: K_0\longrightarrow \mathbb{R}^{2d}$ be defined by
$\psi_t(i)=\psi(i)$ for every $i\neq u$ and
$\psi_t(u)=\psi(v)+t(\psi(u)-\psi(v))$. Thus $\psi_1=\psi$,
and $\lim_{t\mapsto 0}(\sspan(\psi_t(u)-\psi_t(v)))=\sspan(\psi(u)-\psi(v))$. Let
$\psi_0=\lim_{t\mapsto 0}\psi_t$.

Let $\psi_{K,t}^{2d}: \oplus_{T\in K_{d-1}}\mathbb{R}T
\longrightarrow \oplus_{F\in
\binom{K_0}{d-1}}\mathbb{R}^{2d}/\sspan(\psi_t(F))$ be the map induced
by $\psi_t$; thus $\psi_{K,1}^{2d}=\psi_{K}^{2d}$. Let
$\psi_0^{2d}$ be the limit map $\lim_{t\mapsto 0}\psi_{K,t}^{2d}$. Thus for $T$ such that
$\{u,v\}\subseteq T\in K_{d-1}$,
$$\psi_0^{2d}(T)|_{T\setminus
v}=(\psi(u)-\psi(v)) + \sspan(\psi(T\setminus
u))=-\psi_0^{2d}(T)|_{T\setminus u}.$$
Assume for a moment
that $\psi_0^{2d}$ is injective. Then for a small enough perturbation of the
entries of a representing matrix of $\psi_0^{2d}$, the columns of the resulted
matrix would be independent, i.e. the corresponding linear
transformation would be injective. In particular, there would exist
an $\epsilon>0$ such that for every $0<t<\epsilon$, $\Ker
\psi_{K,t}^{2d}=0$, and hence for every $\phi: K_0\longrightarrow
\mathbb{R}^{2d}$ in some Zariski non-empty open set of maps,
$\Ker \phi_K^{2d}=0$. Thus, the following Lemma \ref{lemCol(A)ind} completes the proof.
$\square$

\begin{lemma}\label{lemCol(A)ind}
$\psi_0^{2d}$ is injective for a non-empty Zariski open set of
maps $\psi: K_0\longrightarrow \mathbb{R}^{2d}$.
\end{lemma}
$Proof:$ For every $0<t\leq 1$ and every $F$ such that
$\{u,v\}\subseteq F \in \binom{K_0}{d-1}$,
$\sspan(\psi_t(F))=\sspan(\psi(F))$, and hence in the range of $\psi_0^{2d}$ we
mod out by $\sspan(\psi(F))$ for summands with such $F$. For summands
of $\{u,v\}\nsubseteq F \in \binom{K_0}{d-1}$,
we mod out by $\sspan(\psi_0(F))$. Note that for $T$ such that
$\{u,v\}\subseteq T\in K_{d-1}$,
$\psi_0^{2d}(T)|_{T\setminus
v}=-\psi_0^{2d}(T)|_{T\setminus u}$.

For a linear transformation
$C$, denote by $[C]$ its representing matrix w.r.t. given bases. In
$[\psi_0^{2d}]$ bases are indexed by sets as in
(\ref{eq:HighRigidity}). First add rows $F'\uplus\{u\}$ to rows
$F'\uplus\{v\}$ (in particular $F'\cap \{u,v\}=\emptyset$), then delete the rows $F$ containing $u$, to obtain
a matrix $[B]$, of a linear transformation $B$. In particular, we
delete all rows $F$ such that $\{u,v\}\subseteq F$.

Note that $K'_0=K_0\setminus \{u\}$, thus, for the obvious bases,
$[B]$ is obtained from $[(\psi|_{K'_0})_{K'}^{2d}]$ by doubling the
columns indexed by $T'\uplus\{v\}\in K'_{d-1}$ where both
$T'\uplus\{v\},T'\uplus\{u\}\in K_{d-1}$, and by adding a zero
column for every $T'\uplus\{u,v\}\in K_{d-1}$. For short, denote
$\psi_{K'}^{2d}=(\psi|_{K'_0})_{K'}^{2d}$. More precisely, the
linear maps $B$ and $\psi_{K'}^{2d}$ are related as follows: they
have the same range. The domain of $B$ is $\dom(B)=\dom(\psi_0^{2d})=D_1\oplus
D_2\oplus D_3$ where
\newline $D_1=\oplus\{\mathbb{R}T: T\in K_{d-1}, \{u,v\}\nsubseteq T, (u\in T)\Rightarrow (T\setminus u)\cup v \notin K\}$,
\newline $D_2=\oplus\{\mathbb{R}T: T\in K_{d-1}, u\in T, v\notin T, (T\setminus u)\cup v \in K\}$,
\newline $D_3=\oplus\{\mathbb{R}T: T\in K_{d-1}, \{u,v\}\subseteq T\}$.
\newline
For a base element $1T$ of $D_1$, let $T'\in K'$ be obtained from
$T$ by replacing $u$ with $v$. Then $B(1T)=\psi_{K'}^{2d}(1T')$;
thus $\Ker B|_{D_1}\cong \Ker \psi_{K'}^{2d}$. For a base element $1T$ of
$D_2$, $B(1T)=\psi_{K'}^{2d}(1((T\setminus u)\cup v))$, and
$B|_{D_3}=0$.

Assume we have a linear dependency $\sum_{T\in K_{d-1}}\alpha_T
\psi_0^{2d}(T)=0$.
By assumption, $\{d+2,d+3,...,2d+1\}\notin \Delta^{s}(K')$, hence $\Ker
\psi_{K'}^{2d}=0$, thus $\alpha_T=0$ for every base element $T$
except possibly for $T$ containing $\{u,v\}$ and for $T'\uplus \{u\}, T'\uplus
\{v\} \in K_{d-1}$, where $\alpha_{T'\uplus \{u\}}=-\alpha_{T'\uplus
\{v\}}$.

Let $\psi_0^{2d}|_{\res}$ be the restriction of $\psi_0^{2d}$ to the
subspace spanned by the base elements $T$ such that $v\in T$ and for
which it is (yet) not known that $\alpha_T=0$, followed by
projection into the subspace spanned by the $F\in \binom{K_0}{d-1}$
coordinates where $v\in F$ (just forget the other coordinates). As
$\psi_0^{2d}(T)|_F=0$ whenever $F \ni v \notin T$, if
$\psi_0^{2d}|_{\res}$ is injective, then $\alpha_T=0$ for all $T \in
K_{d-1}$. Thus, the Lemma \ref{lemCol(A|_res)ind} below completes
the proof. $\square$

\begin{lemma}\label{lemCol(A|_res)ind}
$\psi_0^{2d}|_{\res}$ is injective for a non-empty Zariski open set of
maps $\psi: K_0\longrightarrow \mathbb{R}^{2d}$.
\end{lemma}
$Proof:$ Let $G=(\{u\}*(\lk(u,K)\cap \lk(v,K)))_{\leq d-2}$. Note that
$v$ appears in the index set of every row and every column of
$[\psi_0^{2d}|_{\res}]$. Omitting $v$ from the indices of both of the bases
used to define $\psi_0^{2d}|_{\res}$, we notice that
$$\psi_0^{2d}|_{\res}\cong\overline{\psi_0^{2d}|_{\res}}: \oplus_{T\in G_{d-2}}\mathbb{R}T \longrightarrow
\oplus_{F\in
\binom{G_0}{d-2}}\mathbb{R}^{2d}/\sspan(\psi(F\uplus\{v\}))=$$
$$\oplus_{F\in \binom{G_0}{d-2}}(\mathbb{R}^{2d}/\sspan(\psi(v)))/\overline{span(\psi(F))},$$
$$1T\mapsto \sum_{F\in \binom{G_0}{d-2}}\delta_{F\subseteq T} \overline{\psi(T\setminus F)}F$$
where $\delta_{F\subseteq T}$ equals $1$ if $F\subseteq T$ and $0$
otherwise, and $\overline{\sspan(\psi(F))}$ is the image of
$\sspan(\psi(F))$ in the quotient space
$\mathbb{R}^{2d}/\sspan(\psi(v))$.

Consider the projection $\pi: \mathbb{R}^{2d} \longrightarrow
\mathbb{R}^{2d}/\sspan(\psi(v)) \cong \mathbb{R}^{2d-1}$. Let
$\bar{\psi}=\pi \circ \psi|_{G_0}: G_0 \longrightarrow
\mathbb{R}^{2d-1}$, and $\bar{\psi}_G^{2d-1}$ be the induced map
as defined in (\ref{eq:HighRigidity}).
Then $\pi$ induces $\pi_*\overline{\psi_0^{2d}|_{\res}} = \bar{\psi}_G^{2d-1}$.

By assumption, $\{d+1,...,2d-1\}\notin \Delta^{s}(\lk(u,K)\cap
\lk(v,K))$. As symmetric shifting commutes with constructing a cone
(Kalai \cite[Theorem 2.2.8]{skira}, and Babson, Novik and Thomas
\cite[Theorem 3.7]{Babson-Novik-Thomas-Cone}), $\{d+2,...,2d\} \notin
\Delta^{s}(G)$. Hence $y_{2d}^{d-1} \notin \GIN(G)$,
and by Lee \cite{Lee}, $\Ker \phi_G^{2d-1}=0$ for a generic $\phi$.
Thus, all liftings $\psi: K_0\longrightarrow \mathbb{R}^{2d}$ such
that $\bar{\psi}=\phi$ satisfy $\Ker \psi_0^{2d}|_{\res}\cong \Ker
\phi_G^{2d-1}=0$, and this set of liftings is a non-empty Zariski
open set. $\square$

Clearly the set of all $\psi$ such that
$\psi_K^{2d}$ is injective is Zariski open. We exhibited conditions
under which it is non-empty. The choice $k=\mathbb{R}$ was needed for the perturbation argument.
\\
\\
\emph{Proof for exterior shifting}: The proof is similar to the
proof for the symmetric case. We indicate the differences.
$\psi:K_0\rightarrow \mathbb{R}^{d+1}$ defines the first $d+1$
generic $f_i$'s w.r.t. the $e_i$'s basis of $\mathbb{R}^{|K_0|}$ and
induces the following map:
\begin{equation}\label{eq:extHighRigidity}
\psi^{d+1}_{K,\ext}: \oplus_{T\in K_{d-1}}\mathbb{R}T \longrightarrow
\oplus_{1\leq i\leq d+1}\oplus_{F\in \binom{K_0}{d-1}}\mathbb{R}F,\ m\mapsto(f_1\lfloor m,...,f_{d+1}\lfloor m)
\end{equation}
where $f_i\lfloor\cdot$ is the left interior product given by bilinear extension of $e_S\lfloor e_T=\delta_{S\subseteq T}\sgn(S,T)e_{T\setminus S}$, as in \cite{56}.
By \cite[Proposition 3.1]{Nevo-ABC}, $\Ker \psi^{d+1}_{K,\ext}=\cap_{1\leq
i\leq d+1}\Ker f_i\lfloor =
\cap_{R<_{lex}\{d+2,...,2d+1\}}\Ker(f_R\lfloor :\ \oplus_{T\in K_{d-1}}\mathbb{R}T \longrightarrow \mathbb{R})$, and hence by
shiftedness $\{d+2,...,2d+1\}\notin \Delta^e(K) \Leftrightarrow
\Ker \psi^{d+1}_{K,\ext}=0$.

Replacing $\psi(u)$ by $\psi(v)$ induces a map
$$\psi^{d+1}_{K,u}:\oplus_{T\in K_{d-1}}\mathbb{R}T \longrightarrow
\oplus_{1\leq i\leq d+1}\oplus_{F\in \binom{K_0}{d-1}}\mathbb{R}F.$$
By perturbation, if $\Ker \psi^{d+1}_{K,u}=0$ then $\Ker
\psi^{d+1}_{K,\ext}=0$ for generic $\psi$.

Let $[B_{\ext}]$ be obtained from the matrix $[\psi^{d+1}_{K,u}]$ by
adding the rows $F'\uplus u$ to the corresponding rows $F'\uplus v$
and deleting the rows $F$ with $\{u,v\}\subseteq F$. The domain of
$B_{\ext}$ is $D_1\oplus D_2\oplus D_3$ defined by sets indexing a basis as for $B$ in the symmetric
case. For a base element $1T$ of $D_1$, let $T'\in K'$ be obtained
from $T$ by replacing $u$ with $v$. Then
$B_{\ext}(1T)=\psi_{K',\ext}^{d+1}(1T')$; thus $\Ker
B_{\ext}|_{D_1}\cong \Ker \psi_{K',\ext}^{d+1}$. For a base element
$1T$ of $D_2$, $B_{\ext}(1T)=\psi_{K',\ext}^{d+1}(1((T\setminus
u)\cup v))$, and as we may number $v=1,u=2$ then $B|_{D_3}=0$ (the
rows of $F'\uplus u$ and of $F'\uplus v$ have opposite sign in
$\psi^{d+1}_{K,u}$). Now we can adopt the arguments showing that
$\Ker \psi^{2d}_0=0$ using $B$ in the symmetric case, to
show that $\Ker \psi^{d+1}_{K,u}=0$ using $B_{\ext}$.
$\square$

\begin{corollary}\label{cor:contraction}
Let $K$ be a $2d$-sphere for some $d\geq 1$, and let $a,b\in K$ be
two vertices which satisfy the \emph{Link Condition}, i.e that
$\rm{lk}(a,K)\cap\rm{lk}(b,K)=\rm{lk}(\{a,b\},K)$.
Let $K'$ be obtained from $K$ by contracting $a\mapsto b$.
Then:

(1) $K'$ is a $2d$-sphere, PL homeomorphic to $K$ (\cite[Theorem 1.4]{Nevo-VK}).

(2) If $K'$ is $d$-Lefschetz and $\lk(\{a,b\},K)$ is
$(d-1)$-Lefschetz over $\mathbb{R}$, then $K$ is $d$-Lefschetz over $\mathbb{R}$ (by Proposition \ref{propAlgContruction}).
$\square$
\end{corollary}

Let $K$ be a simplicial complex. Its \emph{Stellar subdivision at a
face $T\in K$} is the operation $K\mapsto K'$ where
$K'=\Stellar(T,K):=(K\setminus \st(T,K))\cup(\{v_T\}*\partial T*
\lk(T,K))$, where $v_T$ is a vertex not in $K$
and $\st(T,K)=\{S\in K: T\subseteq S\}$.
Note that for $u\in
T\in K$, $u,v_T\in K'$ satisfy the Link Condition and their
identification results in $K$. Further,
$\lk(\{u,v_T\},K')=\lk(u,\partial T* \lk(T,K))=\partial(T\setminus
\{u\})* \lk(T,K)$.
\\
\\
\emph{Proof of Theorem \ref{thm:connsum,join,Stellar}(3):} Let $T=\Stellar(F,K)$, denote its dimension by $d-1$, and
assume by contradiction that $T$ is not SL. As we have seen in the
proof of Lemma \ref{lem:WL->SL}, there exists $0\leq i_0\leq \lfloor
d/2\rfloor$ such that $\Omega_{SL}(T,i_0)=\emptyset$. First we show
that $i_0\neq \lfloor d/2\rfloor$: for $d$ even this is obvious. For
$d$ odd, note that for $u\in F$ the contraction $v_F\mapsto u$ in
$T$ results in $K$, which is $\lfloor d/2\rfloor$-Lefschetz.
Further, the $(d-3)$-sphere
$\lk(\{v_F,u\},T)=\lk(F,K)*\partial(F\setminus \{u\})$ is SL by
Theorem \ref{thm:connsum,join,Stellar}(1), and in particular is $(\lfloor
d/2\rfloor-1)$-Lefschetz. Thus, by Corollary \ref{cor:contraction},
$T$ is $\lfloor d/2\rfloor$-Lefschetz, and hence $0\leq i_0\leq
\lfloor d/2\rfloor-1$.

Let $L=T*\partial \sigma$, where $\sigma$ is the
$(d-2i_0-1)$-simplex (then $L$ has even dimension $2d-2i_0-2$). By
Lemma \ref{lem:WL->SL}, for any two l.s.o.p.'s $\Theta_T$ and
$\Theta_{\partial \sigma}$ of $\mathbb{R}[T]$ and
$\mathbb{R}[\partial \sigma]$ respectively, $\mathbb{R}[L]/(\Theta_T
\cup \Theta_{\partial \sigma})$ has no $(d-i_0-1)$-Lefschetz
element.

On the other hand, we shall now prove the existence of such
l.s.o.p.'s and a $(d-i_0-1)$-Lefschetz element, to reach a
contradiction. This requires a close look on the proof of
Proposition \ref{propAlgContruction}.

Note that $L=\Stellar(F,K*\partial \sigma)$, and that for $u\in F$
the contraction $v_F\mapsto u$ in $L$ results in $K*\partial
\sigma$. Further, $\lk(\{v_F,u\},L)=\lk(F,K)*\partial(F\setminus
\{u\})* \partial \sigma$.

Applying Zariski topology considerations to subspaces of the space
of maps $\{f:L_0\rightarrow \mathbb{R}^{2d-2i_0}\}\cong
\mathbb{R}^{|L_0|\times (2d-2i_0)}$, we now show that there exists
a map $\psi: L_0\longrightarrow \mathbb{R}^{d}\oplus
\mathbb{R}^{d-2i_0-1}\oplus \mathbb{R}$ such that the following
three properties hold \emph{simultaneously}:

(1) $\psi(K_0) \subseteq \mathbb{R}^{d}\oplus 0 \oplus \mathbb{R}$
and induces an l.s.o.p. $\Theta_K$ of $\mathbb{R}[K]$ (by first $d$ columns) and an SL
element $\omega_K$ of $\mathbb{R}[K]/(\Theta_K)$ (by last column); $\psi(\sigma_0)
\subseteq 0 \oplus \mathbb{R}^{d-2i_0-1}\oplus \mathbb{R}$ and
induces an l.s.o.p. $\Theta_{\partial \sigma}$ of
$\mathbb{R}[\partial \sigma]$ and an SL element $\omega_{\partial
\sigma}$ of $\mathbb{R}[K]/(\Theta_{\partial \sigma})$ (by last $d-2i_0$ columns). By Theorem
\ref{thm:*}, $\omega_K+\omega_{\partial \sigma}$ is an SL element of
$\mathbb{R}[K*\partial \sigma]/(\Theta_K \cup \Theta_{\partial
\sigma})$.

In matrix language, the first $2d-2i_0-1$ columns of
$[\psi|_{K_0\cup \sigma_0}]$ form an l.s.o.p. of
$\mathbb{R}[K*\partial \sigma]$, and its last column is the
corresponding SL element.

(2) $0\neq \psi(v_F)\in \mathbb{R}^{d}\oplus 0 \oplus \mathbb{R}$
induces a map $\pi: \mathbb{R}^{2d-2i_0}\rightarrow
\mathbb{R}^{2d-2i_0}/\sspan \psi(v_F) \cong \mathbb{R}^{2d-2i_0-1}$
such that $\pi\circ\psi|_{K_0\cup \sigma_0}$ induces an element in
$\Omega(G,d-i_0-2)$ for $G=\{u\}* \lk(\{v_F,u\},L)$.

To see this, consider e.g. a map $\psi'$ with
$\psi'(v_F)=(1,0,...,0)$, $\psi'(u)=(0,1,0,...,0)$, $\psi'(s)$
vanishes on the first two coordinates for any $s\in K_0\setminus
\{u\}$ and in addition $[\psi']$ vanishes on all entries on which we
required in (1) that $[\psi]$ vanishes. By Theorem \ref{thm:*} there
exists such $\psi'$ so that its composition with the projection
$\pi':\mathbb{R}^{2d-2i_0}\rightarrow \mathbb{R}^{2d-2i_0}/\sspan
\{\psi(v_F),\psi(u)\}$ induces a pair $(\Theta,\omega)$ of an
l.s.o.p. and an SL element for
$\lk(\{v_F,u\},L)=\lk(F,K)*\partial(F\setminus \{u\})*
\partial \sigma$. By adding $x_u$ to this l.s.o.p. we obtain an l.s.o.p.
for $G$ where $\omega: H(G)_{d-i_0-2}\rightarrow H(G)_{d-i_0-1}$ is
injective; hence property (2) holds for $\psi'$.

The restriction of maps $\psi$ with property (2) to $\st(F,K)_0\cup
\{v_F\}$ is a nonempty Zariski open set in the space of maps
$\{f: \st(F,K)_0\cup \{v_F\}\rightarrow \mathbb{R}^{d}\oplus 0
\oplus \mathbb{R}\}$. The restriction of maps $\psi$ with property
(1) to $K_0$ is a nonempty Zariski open set in the space of
maps $\{f: K_0\rightarrow \mathbb{R}^{d}\oplus 0 \oplus
\mathbb{R}\}$. Hence, their projections on the linear subspace $\{f:
\st(F,K)_0\rightarrow \mathbb{R}^{d}\oplus 0 \oplus \mathbb{R}\}$
are nonempty Zariski open sets (in this subspace). The intersection
of these projections is again a nonempty Zariski open set, thus
there are maps $\psi$ for which both properties (1) and (2) hold.

(3) $\psi(K_0\cup \{v_F\}) \subseteq \mathbb{R}^{d}\oplus 0 \oplus
\mathbb{R}$ and the first $d$ columns of $[\psi]$ induce an l.s.o.p.
$\Theta_T$ of $\mathbb{R}[T]$.

The set of restrictions $\psi|_{T_0}$ of maps $\psi$ with property
(3) is nonempty Zariski open in the subspace $\{f: T_0\rightarrow
\mathbb{R}^{d}\oplus 0 \oplus 0\}$; hence, so is its projection on
the linear subspace $\{f: \st(F,K)_0\rightarrow \mathbb{R}^{d}\oplus
0 \oplus 0\}$. By similar considerations to the above, there are
maps $\psi$ for which all the properties (1), (2) and (3) hold.

The proof of Proposition \ref{propAlgContruction} together with
properties (1) and (2) tell us that for small enough $\epsilon$, the
map $\psi": L_0\longrightarrow \mathbb{R}^{2d-2i_0}$ defined by
$\psi"(v_F)=\psi(u)+\epsilon (\psi(v_F)-\psi(u))$ and
$\psi"(v)=\psi(v)$ for every other vertex $v\in L_0$, satisfies $\Ker
\psi"_{L}^{2d-2i_0}=0$ (see equation (\ref{eq:HighRigidity}) for the
definition of this map). As a nonempty Zariski open set is dense, by
looking on the subspace of maps $\{f:T_0\rightarrow \mathbb{R}^{d}\oplus 0
\oplus \mathbb{R}\}\}$ , we can take $\psi(v_F)$ and $\epsilon$ such
that $\psi"$ satisfies property (3) as well.

Thus, the first $d$ columns of $[\psi"]$ induce an l.s.o.p.
$\Theta_T$ of $T$, the next $d-i_0-1$ columns induce an l.s.o.p.
$\Theta_{\partial \sigma}$ of $\partial \sigma$, and the last column
of $[\psi"]$ is a $(d-i_0-1)$-Lefschetz element of
$\mathbb{R}[L]/(\Theta_T \cup \Theta_{\partial \sigma})$. This
contradicts our earlier conclusion, which was based on assuming that
the assertion of this theorem is incorrect. $\square$

\begin{corollary}\label{cor:Stellar}
Let $\mathcal{S}$ be a family of homology spheres which is closed
under taking links and such that all of its elements are SL, over $\mathbb{R}$. Let
$\mathbb{S}=\mathbb{S}(\mathcal{S})$ be the family obtained from
$\mathcal{S}\cup \{\partial \sigma^n: n\geq 1\}$ by taking the
closure under the operations: (0) taking links; (1) join; (2)
Stellar subdivisions. Then every element in $\mathbb{S}$ is SL.
\end{corollary}
$Proof:$ We prove by double induction - on dimension, and on the
sequence of operations of types (0),(1) and (2) which define $S\in
\mathbb{S}$ - that $S$ and all its face links are SL. Let us call
$S$ with this property \emph{hereditary SL}.

Note that every $S\in \mathcal{S}$ and every boundary of a simplex, is hereditary SL. This includes the (unique) zero-dimensional sphere and provides the base of the induction. (Actually it is known that every (homology) sphere of dimension $\leq 2$ is hereditary SL.)

Clearly if $S$ is hereditary SL, then so are all of its links, as
$\lk(Q,(\lk(F,S))=\lk(Q\uplus F,S)$. If $S$ and $S'$ are hereditary
SL then by Theorem \ref{thm:*} so is $S*S'$ (here we note that every
$T\in S*S'$ is of the form $T=F\uplus F'$ where $F\in S$ and $F'\in
S'$, and that $\lk(T,S*S')=\lk(F,S)*\lk(F',S')$). We are left to
show that if $F\in S$ and $S$ is hereditary SL, then so is
$T:=\Stellar(F,S)$. Assume $\ddim F\geq 1$, otherwise there is
nothing to prove. First we note that by the induction hypothesis for
every $v\in T_0$, $\lk(v,T)$ is hereditary SL:
\\ Case $v=v_F$: $\lk(v_F,T)=\lk(F,S)*\partial F$ is hereditary SL
by Theorem \ref{thm:*}, as argued above.
\\ Case $v\in F$: $\lk(v,T)=\Stellar(F\setminus \{v\},\lk(v,S))$ is hereditary SL
by the induction hypothesis on the dimension.
\\ Case $v\notin F$, $v\neq v_F$ and $F\in \lk(v,S)$:
$\lk(v,T)=\Stellar(F,\lk(v,S))$ is hereditary SL by the induction
hypothesis on the dimension.
\\ Otherwise: $\lk(v,T)=\lk(v,S)$ is hereditary SL.

We are left to show that $T$ is SL: $S$ is SL, and for $u\in F$
$\lk(\{v_F,u\},T)=\lk(F,S)*\partial(F\setminus \{u\})$ is SL by
Theorem \ref{thm:*}. Thus, by Theorem \ref{thm:connsum,join,Stellar}(3) $T$
is SL, and together with the above, $T$ is hereditary SL. $\square$
\\
 \\
The barycentric subdivision of a simplicial complex
$K$ can be obtained by a sequence of Stellar subdivisions: order the
faces of $K$ of dimension $>0$ by weakly decreasing size, and
perform Stellar subdivisions at those faces according to this order;
the barycentric subdivision of $K$ is obtained. Brenti and Welker
\cite[Corollary 3.5]{Brenti-Welker} showed that the $h$-polynomial
of the barycentric subdivision of a Cohen-Macaulay complex has only
simple and real roots, and hence is unimodal. In particular,
barycentric subdivision preserves non-negativity of the $g$-vector
for spheres with all links being SL. The above corollary shows that
the hereditary SL property itself is preserved.

\section{Lefschetz properties and connected sum}\label{sec:WL,SL&Connected Sum}
Let $K$ and $L$ be pure simplicial complexes which intersect in a
common closed facet $<\sigma>=K\cap L$. Their \emph{connected sum over
$\sigma$} is $K\#_{\sigma}L=(K\cup L)\setminus \{\sigma\}$.

\begin{theorem}\label{thm:connected-sum}
Let $K$ and $L$ be homology $(d-1)$-spheres over a field $\mathbb{F}$ which
intersect in a common closed facet $<\sigma>=K\cap L$.
Let $A=\mathbb{F}[x_v:v\in (K\cup L)_0]$. Then:

(0) $K\#_{\sigma}L$ is a homology $(d-1)$-sphere; in
particular its $h$-vector is symmetric.

(1) Let $\Theta$ be a common l.s.o.p for $K$, $L$, $<\sigma>$ and
$K\#_{\sigma}L$ over $A$ (it exists if $\mathbb{F}$ is infinite).
Assume that $K$ and $L$ are $i$-Lefschetz for some $i>0$
 and let $\omega$ be an
$i$-Lefschetz element for both $K$ and $L$ w.r.t. $\Theta$ (it exists). Then $\omega$ is an
$i$-Lefschetz element of $\mathbb{F}[K\#_{\sigma}L]/(\Theta)$.
\end{theorem}

$Proof$: Straightforward Mayer-Vietoris and Euler characteristic
arguments show that $K\#_{\sigma}L$ is a homology $(d-1)$-sphere.

For a simplicial complex $L$ let $\mathbb{F}(L):=\bigoplus_{a:
\rm{supp} (a)\in L}\mathbb{F}x^a$ be a module over $A_{L_0}=\mathbb{F}[x_v: v\in
L_0]$ defined by $x_v(x^a)=\{^{x_vx^a \ \rm{if}\ v\cup
\rm{supp} (a)\in L}_{0 \ \rm{otherwise}}$.
Note that $\mathbb{F}(L)\cong \mathbb{F}[L]$ as $A_{L_0}$-modules.
For $v\in (K\cup L)_0\setminus L_0$ and $m\in \mathbb{F}(L)$, $x_v m=0$.

Then the following is an exact sequence of $A$-modules:
\begin{equation}\label{eq:conn-sum-exact}
\begin{CD}
0 \rightarrow \mathbb{F}(<\sigma>) @>(\iota,-\iota)>>
(\mathbb{F}(K)\oplus \mathbb{F}(L)) @>\iota_K+\iota_L>>
\mathbb{F}(K\cup_{\sigma}L) \rightarrow 0
\end{CD}
\end{equation}
where the $\iota$'s denote the obvious inclusions.
$|\mathbb{F}|=\infty$ guarantees the existence of an l.s.o.p. for each of the $(d-1)$-complexes in Theorem \ref{thm:connected-sum}(1), and as a finite
intersection of Zariski nonempty open sets is nonempty, $\Theta$ as
in (1) exists. When we mod out
$\Theta$ from (\ref{eq:conn-sum-exact}), which is the same as tensor
(\ref{eq:conn-sum-exact}) with $\otimes_A A/\Theta$, we obtain an
exact sequence of $A$-modules:
\begin{equation}\label{eq:conn-sum-modTheta-exact}
\frac{\mathbb{F}(<\sigma>)}{(\Theta)\mathbb{F}(<\sigma>)}\rightarrow
\frac{\mathbb{F}(K)}{(\Theta)\mathbb{F}(K)}\oplus \frac{\mathbb{F}(L)}{(\Theta)\mathbb{F}(L)} \rightarrow
\frac{\mathbb{F}(K\cup_{\sigma}L)}{(\Theta)\mathbb{F}(K\cup_{\sigma}L)}\rightarrow 0
\end{equation}
where in the middle term we used distributivity of $\otimes$ and
$\oplus$.  Note that $\frac{\mathbb{F}(<\sigma>)}{(\Theta)\mathbb{F}(<\sigma>)}\cong \mathbb{F}$
is concentrated in degree $0$ and that
$(\mathbb{F}(K\#_{\sigma}L)/(\Theta))_{<d} \cong
(\mathbb{F}(K\cup_{\sigma}L)/(\Theta))_{<d}$. Thus, for $0<i\leq
d/2$ we obtain the following commutative diagram of $A$-modules:
\begin{equation}\label{eq:conn-sum-omega-isom}
\begin{CD}
(\frac{\mathbb{F}(K\#_{\sigma}L)}{(\Theta)\mathbb{F}(K\#_{\sigma}L)})_i@>\cong >>(\frac{\mathbb{F}(K\cup_{\sigma}L)}{(\Theta)\mathbb{F}(K\cup_{\sigma}L)})_i@>\cong >>(\frac{\mathbb{F}(K)}{(\Theta)\mathbb{F}(K)})_i\bigoplus (\frac{\mathbb{F}(L)}{(\Theta)\mathbb{F}(L)})_{i}\\
@VV\omega^{d-2i} V @VV\omega^{d-2i} V @VV\omega^{d-2i}\oplus \omega^{d-2i}V \\
(\frac{\mathbb{F}(K\#_{\sigma}L)}{(\Theta)\mathbb{F}(K\#_{\sigma}L)})_{d-i}@>\cong >>(\frac{\mathbb{F}(K\cup_{\sigma}L)}{(\Theta)\mathbb{F}(K\cup_{\sigma}L)})_{d-i}@>\cong >>(\frac{\mathbb{F}(K)}{(\Theta)\mathbb{F}(K)})_{d-i}\bigoplus (\frac{\mathbb{F}(L)}{(\Theta)\mathbb{F}(L)})_{d-i} \\
\end{CD}
\end{equation}
where the right vertical arrow is an isomorphism by assumption.
Hence, the left vertical arrow is an isomorphism as well, meaning
that $\omega$ is an $i$-Lefschetz element of
$\mathbb{F}[K\#_{\sigma}L]/(\Theta)$.
$\square$
\\
\\
\emph{Proof of Theorem \ref{thm:connsum,join,Stellar}(2):} If $K$ and $L$ are SL homology $(d-1)$-spheres then by Theorem \ref{thm:connected-sum} $K\# L$ is a homology $(d-1)$-sphere and has a pair $(\Theta,\omega)$ of l.s.o.p. and $i$-Lefschetz element for every $0<i\leq \lfloor d/2\rfloor$.

For $i=0$, as $K\# L$ is Cohen-Macaulay with l.s.o.p. $\Theta$ and
$h_d=1$, then there exists a  $0$-Lefschetz element $\tilde{\omega}$
(i.e. $\tilde{\omega}^{d}\neq 0$. This is equivalent to
$[2,d+1]\in\Delta^s(K\# L)$, which reflects the fact that $K\# L$
has non-vanishing top homology.). By Lemma \ref{lemOmegaZariski} the
sets of $0$-Lefschetz elements and of $(0<)$-Lefschetz elements are
Zariski open. The fact that they are nonempty implies that so is
their intersection, i.e. $K\# L$ is SL.
Similarly, one concludes that if $K$ and $L$ are weak-Lefschetz then so is
$K\# L$.
$\square$
\\
\begin{remark}: The assertion of Theorem \ref{thm:connsum,join,Stellar}(2), rephrased in terms of algebraic shifting, says that if $\Delta^s(K),\Delta^s(L)\subseteq \Delta(d)$ then also $\Delta^s(K\# L)\subseteq \Delta(d)$. The analogues statement for exterior shifting is also true. These assertions follow from the characterization of the algebraic shifting of a union of complexes whose intersection is a simplex, given in \cite{Nevo-ABC}. To obtain the shifting of $K\# L$ from the shifting of $K\cup L$ just delete the facet $\{2,3,...,d,d+2\}$ which represent the extra top homology in $K\cup L$.
\end{remark}
\begin{acknowledgement}
We deeply thank Satoshi Murai for his helpful comments on an earlier version of this paper.
\end{acknowledgement}
\bibliographystyle{amsplain}
\bibliography{biblio,ubt,topology}
\end{document}